\documentclass[12pt]{article}
\usepackage
[pdftex
,a4paper
,hyperindex
=
true
,hyperfigures
=
true
,pagebackref
=
true
,bookmarks
=
true
,pdftitle
=
{TITLE
},pdfsubject
=
{},pdfauthor
=
{AUTHOR
},pdfkeywords={KEYWORDS},pdfpagemode=UseOutlines,menubordercolor=1 1 1,colorlinks=true,urlcolor=blue,citecolor=red,linkcolor=blue]{hyperref} 

\usepackage{yfonts}
\usepackage{amsfonts,amsmath,amssymb,epsfig}
\usepackage{amsthm}

\newcommand{\al}{\alpha}

\newcommand{\cli}{c^{(i)}_{\la}}
\newcommand{\BE}{\begin{equation}}
\newcommand{\EE}{\end{equation}}
\newcommand{\la}{\lambda}
\newcommand{\be}{\beta}
\newcommand{\ga}{\gamma}

\newcommand{\RR}{{\mathbb R}}
\newcommand{\R}{{\mathbb R}}

\newcommand{\ZZ}{{\mathbb Z}}
\newcommand{\NN}{{\mathbb N}}
\newcommand{\N}{{\mathbb N}}
\newcommand{\PP}{{\mathbb P}}

\newcommand{\ep}{\varepsilon}

\newcommand{\La}{\Lambda}

\newcommand{\cjk}{c^{(i)}_{j,k}}

\newcommand{\cla}{c^{(i)}_{\lambda}}
\newcommand{\pla}{\psi^{(i)}_{\lambda}}

\newcommand{\pjk}{\psi^{(i)}_{j,k}}
\newcommand{\ptjk}{\tilde{\psi}^{(i)}_{j,k}}

\newtheorem{lem}{Lemma}[section]
\newtheorem{coro}{Corollary}[section]
\newtheorem{Theo}{Theorem}[section]
\newtheorem{prop}{Proposition}[section]
\newtheorem{defi}{Definition}[section]
\theoremstyle{definition}
\newtheorem{rem}{Remark}[section]
\newcommand{\BP}{\begin{prop}}
\newcommand{\EP}{\end{prop}}
\newcommand{\BC}{\begin{coro}}
\newcommand{\EC}{\end{coro}}
\newcommand{\BL}{\begin{lem}}
\newcommand{\EL}{\end{lem}}
\newcommand{\BD}{\begin{defi}}
\newcommand{\ED}{\end{defi}}
\newcommand{\BT}{\begin{Theo}}
\newcommand{\ET}{\end{Theo}}
\newcommand{\BR}{\begin{rem}}
\newcommand{\ER}{\end{rem}}
\newcommand{\mcjk}{\left|\cjk  \right| }
\newcommand{\mcla}{\left|\cla  \right| }
\newcommand{\mpjk}{\left| \pjk  (x) \right| }
\newcommand{\mcpjk}{\left|\cjk \pjk  (x) \right| }

\newsavebox{\fmbox}
\newenvironment{fmpage}[1]
 {\begin{lrbox}{\fmbox}\begin{minipage}{#1}}
 {\end{minipage}\end{lrbox}\fbox{\usebox{\fmbox}}}

\def\E{{\hbox{I\kern-.2em\hbox{E}}}}

\graphicspath{{Figures/}}

\author{C\'eline Esser\footnote{Corresponding author: Celine.Esser@math.univ-lille1.fr Universit\'e Lille 1, Laboratoire Paul Painlev\'e, Cit\'e Scientifique, 59655 Villeneuve d'Ascq, France. This author is supported by the the Labex CEMPI (ANR-11-LABX-0007-01)} \ and St\'ephane Jaffard\footnote{jaffard@u-pec.fr Universit\'e Paris-Est -  LAMA, UMR 8050, CNRS, UPEC,  61 Av. du Gal. de Gaulle,  F-94010, Cr\'{e}teil, France}}

 \date{\today} 
\title{  Divergence of wavelet series:  \\ A multifractal analysis }

\begin{document}
 \renewcommand{\proofname}{\bf Proof}
 \maketitle

\begin{abstract}
We show the  relevance of a multifractal-type analysis for pointwise convergence and divergence properties of  wavelet series:  Depending on the sequence space which the wavelet coefficients sequence belongs to,  we obtain deterministic upper bounds for the Hausdorff dimensions of the sets of points where a given rate of divergence occurs, and we show that these bounds are generically optimal, according to  several notions of genericity.
\end{abstract}

\tableofcontents


\section{Introduction }

\label{intro}

Pointwise convergence properties of Fourier series  have been a major challenge in analysis, culminating in the  famous Carleson-Hunt theorem.  Later, one direction of research  has been to estimate the ``divergence rate''  of partial sums at  exceptional points where  divergence occurs. The relevant tool  to measure the size of these exceptional sets  is the Hausdorff dimension, thus leading to a {\em multifractal analysis}  of  divergence:
 Denote by $S_n f $ the partial sums of the Fourier series of  a $1$ periodic function $f$, i.e.
 \[ S_n f (x) = \sum_{k=-n}^n c_k \, e^{2i\pi k x} \qquad \mbox{ where} \qquad c_k =  \int_{0}^1 f(t) \, e^{-2i \pi k t} dt , \]
and consider the sets
  $$ E^\beta_f = \left\{ x: \ \limsup_{n \rightarrow \infty} \; n^{-\beta} | S_n f(x) | >0 \right\} . $$
   J.-M. Aubry proved  that, 
 if $  f\in L^p ([0,1])$, and $ \beta >0$,  then $
   \dim (E^\beta_f ) \leq 1-\beta p $ (where $\dim (A) $ denotes the Hausdorff dimension of the set $A$)
  and he showed the   optimality of this result, see  \cite{Aubry06}. This was later extended and refined by F. Bayart and Y. Heurteaux, who, in particular,  showed that optimality holds for generic functions of $L^p$ (in the sense supplied by Baire categories and prevalence), see \cite{BH12,BH14}.

 Such properties were also studied in the setting of wavelet series:    J.-M. Aubry obtained  upper bounds on the dimensions of the sets of points where a given divergence rate  of the wavelet series of an $L^p$ function occurs, see \cite{Aubry06}. Additionally, he showed their  optimality  by a specific construction in the case of the Haar wavelet. This last result was recently extended by F. Bayart and Y. Heurteaux who showed that it holds generically in $L^2$ (in the sense supplied by Baire categories), as a consequence of 
  a general framework that they introduced in order to study multifractal phenomena,  see \cite{BH16}.

 The present  paper   extends these results in several directions: The Hausdorff dimensions of the divergence sets  of wavelet series will be bounded under the assumption that the wavelet coefficients satisfy a  Besov type condition, and the optimality of these bounds will be proved in several  senses: We will show that, generically (in the sense of prevalence, Baire categories and  lineability),  they turn out to be equalities, under general conditions on the wavelet system used. 
 We will also prove irregularity everywhere, which will yield the complete {\em divergence spectrum} of such series.

A {\em wavelet basis}   in $d$ variables  denotes a collection of $2^d -1$  functions $\psi^{(i)}  \in L^2 (\RR^d)$, $i \in \{1, \dots, 2^d-1\}$, such that the collection of functions 
  \[ 2^{dj/2} \psi^{(i)}  (2^j \cdot -k) \quad \mbox{ for }i  \in \{ 1, \dots, 2^d-1\},  \;  j \in \ZZ, \; k \in \ZZ^d  \] forms an orthonormal basis of $L^2 (\RR)$. 
Let 
 \[Ê\pjk (x) =  \psi^{(i)}  (2^j x -k) . \]  
   Convergence properties   can either be studied for  the {\em wavelet expansion}   of a function (or a distribution) $f$  i.e.  for partial sums    of the series  
\BE \label{wavexp} \sum_{j\in \ZZ}   \sum_{i=1}^{2^d-1} \sum_{k\in \ZZ^d} 2^{dj} \langle f |  \pjk \rangle  \pjk (x)   \EE
(note that, if  wavelets have sufficient regularity and decay properties, 
 the duality products $\langle f |  \pjk \rangle $ are well defined even if $f$ is not  a function). 
They can also concern the {\em partial reconstruction operators} 
\BE \label{wavexp44} P_{J,f}(x) = \sum_{j<J} \sum_{i=1}^{2^d-1}  \sum_{k\in \ZZ^d} 2^{dj} \langle f |  \pjk \rangle  \pjk (x)  , \EE 
which can also be rewritten under the form
\[ÊP_{J,f}(x) =  \sum_{k\in \ZZ^d} 2^{dJ} \langle f |  \varphi_{J,k} \rangle  \varphi_{J,k} (x) \qquad \mbox{ where} \quad Ê \varphi_{J,k} (x) =  \varphi  (2^J x -k) \]
for a specific function $\varphi$; this sequence is  referred to as the {\em multiresolution expansion} of $f$. 


 Wavelet expansions have many remarkable properties:  Under mild hypotheses on the
 $\psi^{(i)}$, wavelet bases  are known to form unconditional bases of many functions spaces, such as $L^p$ spaces for $p \in (1, \infty)$,  and most Sobolev or Besov spaces, see \cite{Meyer90} (for such characterizations, the sum corresponding to    negative $j$s in \eqref{wavexp} may have to be replaced by the term $P_{0,f}$). This implies that the wavelet expansion of an element of the space converges with respect to the  corresponding norm. 
However,   the introduction of  the first wavelet basis was  motivated by a  different question:  
   A. Haar defined the {\em Haar basis} in 1909  in order to give an example of an orthonormal basis  for which the expansion of a continuous function converges uniformly (in contradistinction with the Fourier series case). 
 Not surprisingly, for other wavelet bases,   the  expansion of a  bounded continuous function also converges uniformly, see \cite{Walter95}. 
 
 
  Additionally,  the  expansion of an $L^p$ function ($1 \leq p \leq \infty$) converges almost everywhere, and in particular at its Lebesgue points, as shown by S. E. Kelly, M. A. Kon, and L. A. Raphael in \cite{KKR94}; these results were extended by T. Tao to hard and soft sampling in \cite{Tao96}; related pointwise convergence questions were  also investigated in \cite{KV15,Z00}.  These results, however, left open  the study of the rate of divergence at exceptional points  where divergence occurs. 

The point of view that we will adopt  is  different:  We will not consider directly 
the {\em wavelet expansion}  of a function in a  given function space, but rather
{\em  wavelet series}  where the sequence of wavelet coefficients  satisfies a given  (discrete Besov) convergence property.
 These two points of view often coincide because of the wavelet characterizations of function spaces, which are precisely given by such conditions, but it is not always the case: Indeed, on one hand,  some spaces (such as $L^1$ or $L^\infty$) do not have a characterization bearing on the moduli of their wavelet coefficients, and on  other hand, these characterizations only  hold if the wavelets are smooth enough, an assumption that will not be required. Additionally, we will not need  wavelets  to form an orthonormal basis. More precisely, we  will consider  the following setting.    
 
 \BD  Let  $\psi^{(i)}$, $i  \in \{ 1, \dots, N \} $, denote $N$ bounded functions with fast decay defined on $\RR^d$   (where $N$ needs not be equal to $2^d -1$).
 The  associated  wavelet system   is the  collection  of the 
 \[  \pjk (x)=  \psi^{(i)}  (2^j x -k)  \quad \mbox{ for } \quad  i  \in \{ 1, \dots, N\},  \;  j \geq 0, \; k \in \ZZ^d. \]Ê
 \ED  
 

 Formally, a  {\em wavelet series}  will denote a series of the form 
 \BE \label{wavexp3}   \sum_{j\geq 0 } \sum_{i=1}^N \sum_{k\in \ZZ^d} \cjk\pjk (x) \, . \EE
 
\BR  Comparing with  \eqref{wavexp}, we do not  consider the   ``low frequency part''   corresponding to negative values of $j$, 
  the  convergence properties of which are straightforward.  \ER 
  
  \BR We do not assume orthogonality, or vanishing moments, so that this setting  also covers biorthogonal bases, frames,... 
\ER

The purpose of  this paper is to develop  an analysis of  generic convergence and divergence properties of wavelet series.  
In Section \ref{sec:conv} we  define appropriate  notions of   convergence and divergence rates,  and 
  we bound the Hausdorff dimensions  of the corresponding divergence sets  assuming that the coefficients satisfy a   Besov-type condition. In Section \ref{sec:saturatingprocess}, we study a particular example of random wavelet series which will play a key role:  Theorem \ref{th1prev} asserts  that the upper bounds obtained  are generically reached in the sense of  prevalence (which is a natural extension  to infinite dimensional spaces of the concept  of ``almost everywhere'', see Definition \ref{def:prev}).  
   In Section \ref{sec:generic},  these prevalent properties  are shown to hold for  other notions of genericity: Baire categories and lineability.  
In Section \ref{sec:point}, the divergence  rate  of a generic wavelet series  at a  given  arbitrary point  is determined.


\section{Convergence and divergence of wavelet series} 
\label{sec:conv}

In contradistinction with one-variable trigonometric series, which have a  ``natural'' ordering,   simple convergence of a wavelet series at a point does not have a canonical definition; this problem has been  addressed in  \cite{KKR94} where  definitions fitted to simple convergence were introduced, for which the  wavelet expansion of $ L^p$ functions converges almost everywhere.  
%
%
  Such results raise  the problem of precisely estimating the size of the  divergence set   if $f$ is smoother, 
  or  of  determining how fast the wavelet series diverges, when it does. A first  result in this direction stated that, if  $f \in L^{p,s}(\RR^d) $ for an $s >0$, then the wavelet series of $f$ is absolutely convergent outside of a set of dimension at most $d-sp$ (see Prop.  13.6 of \cite{SJdiv}).  
Let $\beta >0$ and let $F^\beta_f$ denote the set of points where the partial sums of the wavelet series (up to scale $J$) are not bounded by $C  2^{\beta J}$; J.-M. Aubry put into light a similarity  between Fourier series and wavelet series, showing that, if $f \in L^p (\RR) $, then $  \dim (F^\beta_f ) \leq 1-\beta p $, see \cite{Aubry06}. 

From now on, we assume that a wavelet system generated by the $\psi^{(i)} $, $i \in \{1, \dots, N\}$, has been fixed.


%

 
%
%



\BD Let $ \ga \in \RR$ $x \in \RR^d$,  ${ \mathcal C}  = \left\{\cjk \right\}$  be a collection of  coefficients and $\left\{\pjk\right\}$ be a wavelet system;  the corresponding series  \eqref{wavexp3} {diverges at rate at least $\ga$} at $x$ if there exist a constant $C>0$ and a sequence $(i_n,j_n,k_n)$ of indices such that $j_n \to + \infty$ and 
\BE \label{divpoi}   \left| c^{(i_n)}_{j_n, k_n} \psi^{(i_n)}_{j_n, k_n} (x) \right|  \geq C 2^{\ga j_n} ,  \EE
which we will denote by $ { \mathcal C}  \in D^\gamma (x) $. 

The  {divergence exponent} of the wavelet series associated with ${ \mathcal C}$ at $x$ is 
\[ \delta_{ \mathcal C}(x) = \sup \{ \gamma : \ { \mathcal C} \in D^\gamma (x)\} , \] 
and the {$\ga$-divergence sets} are the sets 
\[ {\mathcal E}^\ga_{ \mathcal C} =  \{ x : \ \delta_{ \mathcal C} (x) = \ga \} . \] 

The mapping   
\[{\mathcal D}_{ \mathcal C} : \ga \mapsto   \dim \left(   {\mathcal E}^\ga_{ \mathcal C} \right)   \]
 is  the {wavelet divergence  spectrum} of the sequence ${ \mathcal C}$. 
\ED 

\BR
In this definition, the reference to the wavelet system used is implicit. When needed by the context, if $\Psi$ denotes the wavelet system $\left\{\pjk\right\}$ which is  used, we will write $ \delta^{\Psi}_{ \mathcal C}(x)$ for the corresponding divergence exponent.
\ER

Note that (\ref{divpoi}) implies a divergence of the wavelet series only if $\gamma \geq 0$. On other hand, if $\gamma <0$, the negation of  (\ref{divpoi}) actually expresses a convergence rate for the partial reconstruction operators
\[ P_{J ,{ \mathcal C}} (x) =   \sum_{j =0}^J \sum_{i=1}^N \sum_{k \in \ZZ^d}   \cjk \;  \pjk  (x),  \] 
as illustrated in Proposition \ref{propanal}  below. Though we are mainly interested in convergence and divergence properties of wavelet series, 
 we also give results concerning the sequence $P_{J ,{ \mathcal C}}(x)  $. 
 In order to state them,  we introduce a divergence criterium for   multiresolution sequences. 

\BD Let $\gamma >0$ and $x \in \RR^d$; the sequence $P_{j ,{ \mathcal C}} (x)$ diverges at rate $\gamma $ if 
\[ \exists C >0, \; \exists j_n \rightarrow + \infty : \qquad | P_{j_n ,{ \mathcal C}} (x) | \geq C \; 2^{\gamma j_n} . \]
\ED

\BP \label{propanal} Assume that 
\BE  \label{estimco} \exists a \in \RR , \, \exists C_0 >0, \; \forall i,j,k : \qquad  \mcjk \leq C_0 2^{a j}  . \EE

\begin{itemize} 
\item Let $\ga <0$ and $x \in \RR^d$;  
 if 
$$\label{divpoineg} \exists C_1 >0, \;  \forall i,j,k :    \qquad   \mcpjk \leq C_1 2^{\ga j} ,   $$ 
then $ P_{J ,{ \mathcal C}} (x)  $ has a limit $f_{ \mathcal C} (x)$ when $J \rightarrow + \infty$ and 
\BE   \label{estimapprox} \forall \ga' > \ga,\;  \exists C_2 >0, \;   \forall J \geq 0,  \qquad | f_{ \mathcal C} (x) - P_{J ,{ \mathcal C}} (x) | \leq C_2 2^{\ga' J} . \EE

\item Let $\ga >0$ and $x \in \RR^d$;   if the sequence $P_{j ,{ \mathcal C}} (x)$ diverges at rate $\gamma $, then 
$ \forall \delta < \gamma  $, the corresponding wavelet series diverges at rate at least $\delta$ at $x$. 
\end{itemize}
\EP

\begin{proof} We will use the convention that $C$ is a generic constant which only depends on the wavelet system and its  value may change from one line to another. 

We first prove the first point. For every $j \geq 0$, let
\[ Q_{j ,{ \mathcal C}}  (x) =   \sum_{i=1}^N \sum_{k \in \ZZ^d}     \cjk \;  \pjk  (x).  \] 
 First note that, if the wavelets are compactly supported,  then $Q_{j ,{ \mathcal C}}  (x) $ has a bounded number of terms (which only depends on the size of the support), so that $|Q_{j ,{ \mathcal C}}  (x)| \leq C 2^{\ga j} $ and  the sequence  $|Q_{j ,{ \mathcal C}}  (x)| $ decays geometrically; thus  \eqref{estimapprox}  follows with $\ga' = \ga$. Otherwise,
 let $\ep >0$;  we split the sum into two terms. First, if $| k-2^j x| \leq 2^{\ep j} $ the sum has less that  $3 \cdot 2^{ d \ep j} $  terms, and we conclude as in the  compactly supported case. Second, if  $| k-2^j x|  > 2^{\ep j} $, then using the fast decay of the wavelets 
 \[\forall N ,\, \exists C_N> 0 : \quad  \mpjk \leq \frac{C_N}{(1+ | 2^j x-k | )^N} \leq \frac{C_N}{(1+ | 2^j x-k|  )^{N/2}}2^{-\ep Nj/2} \]
 and \eqref{estimco} implies that 
 \[ \sum_{i,k}  \mcpjk  \leq C 2^{-\ep Nj/2} 2^{a j}  ,\]
 which is bounded by $ C 2^{\ga j}$ if we pick $N$ large enough, hence the first point. 
 
 As regards the second point, we note that, if $P_{j ,{ \mathcal C}} (x)$ diverges at rate $\gamma $, then for every $\delta < \gamma$
\[ \exists j'_n \rightarrow + \infty : \qquad \sum_{i=1}^N \sum_{k \in \ZZ^d} \left| c^{(i)}_{j'_n, k} \psi^{(i)}_{j'_n, k} (x)  \right| \geq 2^{\delta j'_n }.\]
Next, we proceed by contradiction assuming that (\ref{divpoi}) does not hold for a rate $\delta<\gamma$. The end of the proof lies then on arguments similar to those used for the proof of the first point. 

\end{proof}

\BR  Sets of divergence of $P_{j ,{ \mathcal C}} $ and of wavelet series  can differ: Indeed, one easily checks that  a wavelet series can diverge at  a given point   at a  rate $\gamma >0$ and nonetheless $P_{j ,{ \mathcal C}}$ may converge at the same point at a    rate $\beta$. 
\ER


\subsection{Upper bound for the wavelet divergence spectrum} 

\label{sec:rates}

We will start by proving upper bounds on the Hausdorff dimensions on the $\al$-divergence sets ${\mathcal E}^\al_{ \mathcal C}$ when ${\mathcal C}$ belongs to  a  discrete Besov space. 
 In order to state them, 
we need to recall the classical wavelet indexing.   If $ k = (k_1, \ldots, k_d) \in \ZZ^d$ and $j \geq 0$, let $\la = \lambda_{j,k} $ denote the dyadic cube of generation $j$:
$$ \label{dyadcub}   \la =  \lambda_{j,k} \; : = \left[ \frac{k_1}{2^j} , \frac{k_1+1}{2^j}   \right) \times \cdots \times  \left[ \frac{k_d}{2^j} , \frac{k_d+1}{2^j}  \right).$$
We will  index  wavelets and wavelet coefficients by $(i,j,k)$ or by $(i, \la )$, writing indifferently $ \cjk$ or $\cli$. We will use similar notations for wavelets, denoting $\psi^{(i)}_{j,k}$ by $\psi^{(i)}_{\lambda}$.  We define $\La_j$ as the set of indices 
$(i, \la )$ such that $\la  $ is of generation $j$.  
We will use the following sequence spaces.

\BD Let $s\in \RR$ and $p,q \in (0,+\infty]$. A sequence ${\mathcal C} = \left\{ \cjk \right\}$  belongs to $b^{s,q}_p $ if  it satisfies 
\BE \label{Besovdis}
 \left( \sum_{(i, \la ) \in \La_j } \left|  \cjk \; 
2^{(s-\frac{d}{p})j}\right|^p
\right)^{1/p} =\ep_j \;\;\; \mbox{with} \;\;\; (\ep_j)_{j \geq 0} \in \ell^q \EE
\ED

 In this definition we do not write down explicitly the dependency of $b^{s,q}_p $ in the space dimension  $d$. Moreover, we use the standard conventions  if $p$ or $q$ take the value $ + \infty$. In particular, ${\mathcal C}\in C^{s} := b^{s ,\infty}_\infty $ if the sequence $\left\{\cjk\;  2^{s j}\right\}$ belongs to
$l^{\infty}$.

\BR The chosen normalization of the coefficients in this definition of the discrete Besov spaces is not the one usually considered, see \cite{P80}: this is justified by the fact that if the wavelets system is composed of smooth enough functions, and forms an orthonormal basis, a biorthogonal basis or a frame, \eqref{Besovdis} characterizes the function space  $B^{s,q}_p (\RR^d)$, see e.g. \cite{Meyer90}. Therefore,  Proposition~\ref{divwase}  below can be reinterperted as yielding bounds on the wavelet divergence exponents of  functions (or distributions) in these spaces. Additionally, the embeddings between Besov and Sobolev spaces implies that the same results hold for functions in the Sobolev space $L^{p,s}  (\RR^d)$. 
\ER

\BP \label{divwase}  Let $s\in \RR$ and  $p, q  \in (0, +\infty] $.  Let $ {\mathcal C}  \in b^{s,q}_p $;  
then  
\[ \forall x \in \RR^d, \qquad \delta_{\mathcal C} (x) \leq \frac{d}{p} -s, \] 
and
\BE \label{maxdiv}\  \forall \ga \in  
\left[  -s ,  -s + \frac{d}{p}  \right] ,  \qquad
\dim \left( \{ x : \; \delta_{\mathcal C}(x) \geq \gamma \}  \right)  \leq d -sp -\ga p. \EE
\EP


\BR The inequality (\ref{maxdiv}) implies that the wavelet divergence spectrum of a sequence $ {\mathcal C}  \in b^{s,q}_p $ satisfies ${\mathcal D}_{\mathcal C}  (\ga) \leq   d -sp - \ga p$. Note also that, when $p =2$,  the  first statement  is in the same spirit as the  results concerning rates of convergence  proved in \cite{KKR2}. 
\ER

\BR \label{rem_holder} The definition of the divergence exponent is very similar to the opposite of an H\"older exponent. This explains why the upper bound (\ref{maxdiv}) obtained for the wavelet divergence spectrum of any sequence of $b^{s,q}_p$ coincides with  the classical upper bound (if one considers it as a function of $- \ga$) of the H\"older spectrum of any function in the corresponding functional Besov space $B^{s,q}_p$, see \cite{Jaffard00}. However, an important  difference is that divergence properties rely on the exact  values of the wavelets at given points, while pointwise H\"older regularity does not depend on the (smooth enough) wavelet basis used. More precisely, the dyadic covering property introduced in Definition \ref{defpsigrand} will play an crucial role in the construction of sequences that saturate the inequality (\ref{maxdiv}).
\ER


\begin{proof}[\bf Proof of Proposition \ref{divwase}] 
For each $q\in (0, + \infty ]$, the space    $b^{s,q}_p$ is included in $b^{s,\infty}_p$ so that it suffices to prove the results for $b^{s,\infty}_p$ i.e. to assume that  the wavelet sequence   satisfies
\BE \label{carbes}  \exists C >0 \, : \quad  2^{-dj} \sum_{(i, \la )  \in \La_j} \mcla^p \leq  C  2^{-spj} \quad \forall j \geq 0. \EE 
This bound implies that each coefficient $\cla$ with $(i,\la)\in \Lambda_j$ satisfies 
$ \mcla \leq C 2^{(d/p-s)j} , $
and since the wavelets $\psi^{(i)}$ are bounded functions, the first statement follows.

If $\la$ is a dyadic cube of width $2^{-j}$, for any $a>0$, let $a \la$  denote the  cube of same center as $\la $, which is homothetical to $\la$ and of  width $a 2^{-j}$. 
Let $\ga \in \RR$.  We define  
\[ E_{j, \ga} = \left\{(i, \la ) \in \Lambda_j : \ \mcla \geq 2^{ \ga j} \right\}  ,  \] 
and, for any $\ep >0$,  
\[ E^\ep_{j, \ga} = \! \bigcup_{\la: \ \exists i : \ (i, \la )  \in E_{j, \ga}} \!\!\! 2^{\ep j } \la   .  \] 
It follows from (\ref{carbes}) that 
\BE \label{cardeja} \text{Card} \: ( E_{j, \gamma}) \leq C  2^{(d-sp-\ga p )j} .  \EE 
Let us set \[ E^\ep_{ \ga}  = \limsup_{j \to + \infty} E^\ep_{j, \ga}. \]
 Because of  (\ref{cardeja}), $E^\ep_{j, \ga}$ is composed of at most $C  2^{(d-sp-\ga p )j}$ cubes of width $2^{-(1-\ep) j} $;
 using these cubes for $j \geq J$ as a  covering of $E^\ep_{ \ga} $, we obtain that 
\BE \label{dimgea}  \dim (E^\ep_{ \ga} ) \leq \frac{d-sp-\ga p }{1-\ep}  . \EE

Let us now estimate $ \left|  \cla   \pla (x) \right| $ for   $x\notin E^\ep_{ \ga}   $. We consider two cases: 
First, assume that   $(i,\la) \notin E_{j, \ga}$; then  $\left|  \cla  \right| < C 2^{\ga j} $, so that 
\BE \label{decal}  |  \cla   \pla (x) |  \leq C' 2^{\ga j} .\EE 

Assume now that   $(i,\la) \in E_{j, \ga}$;  using the fast decay of the wavelets,
\[ \forall N , \; \exists C_N >0 \;  \mbox{ such that} \qquad  | \psi^{(i)} (2^j x-k) | \leq \frac{C_N}{(1 + | 2^j x-k|)^N}  . \] 
Since $x\notin E^\ep_{ \ga}   $ and  $(i,\la) \in E_{j, \ga}$, it follows that, for $j$  large enough,   $| 2^j x-k|  \geq 2^{\ep j}$  where $k \in \ZZ^d$ is such that $\lambda = \lambda_{j,k}$. We obtain that  
\[ | \psi^{(i)} (2^j x-k) | \leq  C_N2^{-\ep N j} . \]
Using  (\ref{carbes}),  $| \cla  | \leq C 2^{-(s-d/p)  j} $,
so that  
\[ | \cla \pla (x) | \leq C_N C 2^{-(s-d/p)  j} 2^{-\ep N j},  \]
 and choosing $N$ large enough yields that (\ref{decal}) also holds.


Thus, outside of $E^\ep_\ga$, the wavelet series  diverges at rate at most $\ga$, i.e. $\delta_{\mathcal C} (x) \leq \gamma$. Since (\ref{dimgea}) holds  for any $\ep >0$, 
\[
\dim (\{x : \delta_{\mathcal C}(x) > \gamma \} )\leq d-sp-\ga p .
\]
This last property holds for any $\gamma \in \R$. Consequently, for a given $\gamma$ and any $\delta>0$, 
\[
\dim(\{x : \delta_f(x) \geq \gamma \} )\leq \dim (\{x : \delta_f(x) > \gamma- \delta \}) \leq d-sp-(\ga - \delta) p 
\]
 and the  upper bounds stated in Proposition \ref{divwase} follows by taking $\delta \to 0$.  
 \end{proof}

 
 \subsection{Maximally divergent wavelet series}
 
A particularly relevant case is supplied by the following  wavelet series,  which are as divergent as allowed by Proposition \ref{divwase}.

\BD \label{defmax} Let $s \in \R$ and $p,q \in (0, + \infty]$. Let $A$ be a  non-empty open subset of $\RR^d$ and let   ${ \mathcal C} \in  b^{s,q}_p $; the corresponding wavelet series (\ref{wavexp3}) is maximally divergent in $A$   if  it  satisfies the following properties: 
\begin{itemize}
\item  for every point $x$ of $A$, $\delta_{ \mathcal C} (x) \in \left[ -s,  \frac{d}{p} -s  \right] ,  $
\item 
for almost every point $x$ of $A$, $\delta_{\mathcal{C}}(x) = -s$, 
 \item for every non-empty open set  $B\subset A$ and for every $\gamma \in \left[ -s,  \frac{d}{p} -s  \right]$, 
 \BE  \label{dimgenerloc}
\dim\left( \{x \in B : \delta_{\mathcal{C}}(x) = \gamma\}\right)  =  d -sp -\ga p  .  \EE
\end{itemize}
\ED

 Note that (\ref{dimgenerloc}) states that the wavelet divergence spectrum is locally invariant  inside $A$ (see \cite{BDJS13} for a precise definition and the basic properties of the related notion of {\em local spectrum}). Remark also that if $p=+ \infty$, then this definition boils down to the condition: $\forall x \in A$,  $\delta_{\mathcal{C}}(x) = -s$. 

We will not only prove the existence of maximally divergent wavelet series (thus proving that the bounds in Proposition \ref{divwase} are sharp), but  we will also show that this extremal behavior  is generic, for several notions of genericity. In order to prove such results, we need that a sufficient number of wavelets do not vanish at a point.  This will be insured by the following property.

\BD \label{defpsigrand} A wavelet system $\left\{\pjk\right\}$ satisfies the dyadic covering property of depth $M$  if there exist $C_0 >0$ and a finite collection of triplets $(i_l, j_l, k_l)$, $l \in \{1, \dots , L\}$,
 such that  $\forall l$, $j_l  >0$  and 
$$ \label{tripl}  \forall x \in [0,1)^d ,\;  \exists  l \in \{Ê1, \dots , L\} : \quad \left| \psi^{(i_l)} (2^{j_l}x-k_l) \right| \geq C_0 ;  $$
the depth $M$  denotes  the largest value taken by the $j_l$. 
 \ED 
 
 This definition will be  the key property in order to obtain divergence results  at \emph{every} point. 
 It is clearly satisfied by the Haar system  with $M =1$ (if properly defined  at dyadic points, e.g. taking $\psi = 1_{[0, 1/2)} -  1_{[1/2, 0)}$ or  taking the definition which ensures  that every point is a Lebesgue point), but it is not satisfied by the Schauder system (divergence does not hold at  dyadic numbers, see the end of Section \ref{sec:conv} for implications of this remark).

  
\BD \label{def_dycov} If a  wavelet system satisfies the dyadic covering property, we define $L$ affine mappings $\mu^l$, $l \in \{1, \dots , L\}$, by the condition that $\mu^l$ is the dilation-translation which maps $[0,1)^d$  to the dyadic cube defined by $(j_l, k_l)$. If $\la$ is an arbitrary dyadic cube, we call the collection of the  $\{ \mu^l (\la) \}_{l \in \{1, \dots , L\}}$ the {dyadic covering of} $\la$. 
 \ED

\BR \label{rem:cov} 
A consequence of this definition  is that, for every $x \in \lambda$, there exists an  $l \in \{1, \dots , L\}$ such that $\left| \psi^{(i_l)}_{\mu^l(\la)} (x) \right| \geq C_0$, and the scale of the dyadic cube $\mu^l(\la)$  exceeds the scale of $\la$ by at most $M$. Note however that the cubes $\mu^l (\la)$ need not  be included in $\la$, which will lead to technicalities in the following proofs; this is in  contradistinction  with the case of the Haar system, which  is  remarkable since one can take for the $\mu_l$  a unique mapping which is the identity (hence simpler proofs in \cite{Aubry06,BH16} for the Haar setting). \ER


We now check that this condition is satisfied by the wavelet systems commonly used.  First, recall that, in a Hilbert space $H$, a sequence of vectors $(e_n)_{n \in \NN} $ is a {\em frame}  if there exist $C, C' >0$ such that 
\[ \forall x \in H, \qquad C \| x\|^2  \leq \sum_{n \in \NN}  \langle x| e_n \rangle^2 \leq C' \| x\|^2 . \]
When such is the case, there exists a sequence $(f_n)$ such that 
\BE \label{frame} x =  \sum_{n \in \NN}   \langle x| f_n \rangle e_n . \EE  The  $(f_n)_{n \in \NN}$ constitute  a {\em dual frame}  of the $(e_n)_{n \in \NN}$  (see e.g. \cite{BenLi,Daub,DGM} for  basics on wavelet frames). 

\BD A wavelet frame system is a couple of wavelet systems $\left\{ \pjk\right\} $ and $\left\{ \ptjk\right\} $  such that the $2^{dj/2}\ptjk$ and  the $2^{dj/2}\pjk$ form  dual frames (for $i=1, \cdots, N$, $j \in \ZZ$, $k \in \ZZ^d$).
\ED

This setting covers a wide range of wavelet expansions, since it is fitted to wavelet orthonormal or bi-orthogonal bases, see \cite{CDF}, or finite unions of such bases, or  dyadic wavelet frames such as developed in \cite{BenLi}.

 \BP \label{lempsigrand}  Any wavelet frame system composed of continuous functions satisfies the dyadic covering property. 
 \EP

\begin{proof} Let us first check that, for a given $x_0\in [0,1]^d$, there exists at least one triplet $(i,j,k)$ with $j\geq 0 $  and  such that 
 \BEÊ\label{triplet} \psi^{(i)} (2^{j}x_0-k) \neq 0. \EE
We will prove this result by contradiction. Let $\ep >0$ and $f_{\ep, x_0} $ be a compactly supported $ C^{\infty}  $  function such that 
\begin{itemize}
\item $Supp ( f_{\ep, x_0} ) \subset B(x_0, \ep )$ 
\item $\forall x$ $f_{\ep, x_0} (x) \in [0,1]$ and $f_{\ep, x_0} (x_0 ) =1$. 
\end{itemize}

 The restatement of \eqref{frame}  for $f_{\ep, x_0} $ in the wavelet frame system setting states that 
  \BE \label{wavexp33}   f_{\ep, x_0} = \sum_{i=1}^N \sum_{j\in \ZZ  }  \sum_{k\in \ZZ^d} \cjk\pjk  \, , \EE
 where 
 \[ \cjk = 2^{dj}  \int f_{\ep , x_0}(x)  \tilde{\psi}^{(i)}  (2^j x -k) \; dx, \]
 where the series converges in $L^2$. However, note that 
 \[ \left| \cjk \right| \leq C \ep^d 2^{dj} \]
so that the sum for $j <0$ converges absolutely and uniformly on any compact set and the corresponding sum is bounded by $C \ep^d$. 
 
 We pick $\ep$ such that $C \ep^d < 1/2$, and we now consider the sum for $j \geq 0$. Since $f_{\ep, x_0}$ is continuous and compactly supported, 
 the sum  converges  uniformly on compact sets.  
 Since $f_{\ep, x_0} (x_0)=1$,  and since the whole series \eqref{wavexp33}  converges to $f_{\ep, x_0}$, it follows that 
 \[ \left| \sum_{i=1}^N \sum_{j\geq 0 } \sum_{k\in \ZZ^d} \cjk \;  \pjk (x_0)  \right| \geq 1/2; \]
 if all the $ \psi^{(i)} (2^{j}x_0-k) $ vanished, then this sum would also vanish; hence a contradiction, 
 and (\ref{triplet}) follows.  \\ 


We now prove Proposition \ref{lempsigrand}. Let $x_0$ be given. Because of the continuity of the wavelet, for the triplet $(i,j,k)$ supplied by (\ref{triplet}), there actually exists an open ball $I_{x_0}$, centered at $x_0$ and on which 
\[  \forall y \in I_{x_0}, \qquad | \psi^{(i)} (2^{j}y-k)  | \geq C_{x_0} >0. \] 
Since the collection of these balls covers $[0,1]^d$, we can extract a finite covering, thus yielding the wavelets in the statement of the lemma. 
\end{proof} 

\BR  The same argument easily applies to more general wavelet frames where the dilation parameters are not dyadic but of the form $a^j$ for an $a >1$ and translation parameters of the form $bka^{-j}$ instead of $k2^{-j}$ (see \cite{DGM}),  as long as the dual frame is of the same form (i.e. of wavelet type). 
\ER

We now show the necessity of the dyadic covering property in order to prove maximal divergence of wavelet series.   We consider the Schauder basis, 
where the generating wavelet is the  function 
 $\La $ defined by 
\[ 
\La (x) = \left\{ 
\begin{array}{ll}
  \min (x, 1-x)  & \mbox{ if } \;\; x \in [0, 1] \, , \\
   0  & \mbox{ else. }
  \end{array} 
  \right.
  \]
The  Schauder basis on the interval $[0, 1]$ is composed of the  functions $1$, $x$ and the  
 \[ \La (2^jx-k) , \qquad \mbox{ for } \;  j \geq 0 \;  \mbox{ and } \;  k \in \{ 0, \dots, 2^j-1 \}; \]
but we can also consider the wavelet system composed of the $ \La (2^jx-k) ,$ for  $ j \geq 0 $ and $ k \in \ZZ$.  This system clearly does not satisfy the  dyadic covering property because, for $j \geq 0$,  all wavelets vanish at the integers. It follows that, {\em whatever the coefficients}  $c_{j,k}$ {Ê\em are}, the series 
\[ \sum_{j \geq 0, k \in\ZZ} c_{j,k} \;  \La (2^jx-k) \] converges at dyadic numbers simply because it has a finite number of nonvanishing terms, and the first condition of maximal divergence does not hold ($\delta_{\mathcal C} = -\infty$ at dyadic points). 

Note that  the same idea  also yields counterexamples in  several variables:  One picks (smooth) wavelets supported in $[0,1]^d$. The corresponding  wavelet series boils down to finite sums on  the  hyperplanes  of equation $x_i = k/2^j$  so that $\delta_{\mathcal C} = -\infty$  on a set of dimension (at least) $d-1$, and the first point of Definition \ref{defmax} never holds. 


\section{Saturating random series and prevalence} 
\label{sec:saturatingprocess}

Up to now, prevalence-type results in multifractal analysis were obtained by the rather technical method of constructing high dimensional {\em probe spaces}, see e.g. \cite{FJ06,Fraysse1}.  Our strategy  will be   more direct, through a  probabilistic construction:  
We  introduce stochastic processes that  allow to obtain prevalent properties about the divergence of wavelet series associated to sequences of $b^{s,q }_p$. More precisely, we construct random wavelet series, referred as {\em saturating series}, whose  coefficients belong to  $b^{s,q }_p$  and which are almost surely maximally divergent. 
We then prove the optimality of Proposition \ref{divwase} in a strong sense: We show how  saturating series yield the prevalent  set of sequences having a maximally divergent wavelet series in  $b^{s,q }_p$.

\subsection{Saturating random series}
Let us assume that a wavelet system satisfying the dyadic covering property of depth $M$ has been chosen. 
From now on, we work on  $[0,1)^d$; we will show afterwards how to extend the results to  $\RR^d$. We fix $s\in \RR$ and $p,q \in (0, + \infty]$. Let $j\geq 1 $ and $k\in \{ 0, \dots , 2^{j}-1 \}^d$ be given, and  denote by  $\lambda$ the corresponding dyadic cube. The integer  $J \leq j$ is defined as follows:  Consider the irreducible representation
\BE \label{irred} \frac{k}{2^j} = \frac{k'}{2^J} \;\;\; \mbox{ where} \;\;\;
k'=(k'_1, \dots k'_d) \EE
and $k'_1, \dots k'_d$ are not all even, and set
\BE \label{defCO1} e^{(i)}_{\la}  = 2^{-(\log j)^2} 2^{(\frac{d}{p}-s)j} 2^{-\frac{d}{p}J}
\EE
for every $i \in \{0, \dots,N\}$. 

\BR \label{hierar} An important remark is that the sequence $\left\{e^{(i)}_{\la} \right\}$ is {\em hierarchical}, i.e.  there exists $\beta \in \RR$ such that $2^{\be j} e^{(i)}_{\la}$  decreases on the dyadic tree when $j$ increases.
\ER

Additionally, let $\xi^{(i)}_{j,k}$ be  I.I.D.  random variables of density 
\BE \label{denslin} \rho(x) = \frac{1_{[-1,1]}(x)}{2} \; dx . \EE
 The random sequence of  coefficients $\cjk$  is  defined as follows: Let  $l >0$ and $\la \subset [0,1)^d$ be a dyadic cube of scale $j  \in \{ m  M +1 , \dots ,(m+1) M\}$.  If there exists  at least one  cube $\nu  \subseteq [0,1)^d$ of scale $mM$ and one $l \in \{1, \dots, L\}$ such that $\mu^l (\nu)= \la$, then we set
 \[  f^{(i)}_\la = \sup \left\{ e^{(i)}_{\nu} : \nu  \subseteq [0,1)^d \mbox{ of scale $mM$ such that  } \exists l \in \{1, \dots, L\} \mbox{ with } \mu^l (\nu)= \la\right\}; \]
otherwise, we set $ f_\la =0$. 
 The purpose of this definition is to ensure that, for  each dyadic cube of $\nu $ of scale $mM$,   the dyadic covering of $\nu$ has coefficients at least equal to $e_\nu$ (note that some of these coefficients  can be larger because a given cube can be the image of different cubes through several mappings), i.e.
\BE \label{majCo}
f^{(i)}_{\mu^l(\nu)} \geq e^{(i)}_\nu
\EE
 for all  $l \in \{1, \dots, L\}$. As mentioned in Remark \ref{rem:cov}, the cube $\mu^l(\nu)$ is not necessary included in $\nu$ (therefore, we do not have a priori relations between $e_{\mu^l(\nu)}$ and $e_\nu$).

 We then define
\BE \label{defCO} \cla = \xi^{(i)}_{\la}f^{(i)}_{\la}
\EE
and we denote by ${\mathcal C}$ the random sequence  $\left\{ \cla\right\}$. 
The aim of this section is to obtain the following result.

\BT \label{thm1}
The random  sequence ${\mathcal C}$ takes values in a compact subset of $b^{s,q }_p$,  and  its associated wavelet series (\ref{wavexp3})  is almost surely  maximally divergent  on $[0,1)^d$.
\ET 

This statement motivates  the name of {\em saturating sequence}, since this random sequence  ``saturates'' the upper bound  of the wavelet divergence spectrum. Theorem~\ref{thm1} will follow from several intermediate results.

\BL \label{regsatgauss}  The random  sequence  ${\mathcal C}$ belongs  to  $b^{s,q}_p$. 
\EL

\begin{proof}  Let us first check that it is the case for the sequence $\left\{e^{(i)}_\la \right\}$.  Let us fix $j \ge 0$. For every $J \le j$, there are less than $C  2^{dJ}$  coefficients satisfying (\ref{irred}); thus, 
\[ \sum_{(i,\la) \in \La_j}  \left| e^{(i)}_\la 2^{(s-\frac{d}{p})j}\right|^p \leq  C 2^{-p (\log j)^2} \sum_{J=0}^j 2^{dJ}\left(2^{-dJ/p}\right)^p= C   j2^{-p (\log j)^2},\]
and \eqref{Besovdis} is satisfied since 
$$ \sum_{j \ge 0} j^{\frac{q}{p}}2^{-q (\log j)^2} <\infty .$$

The construction of the sequence $\left\{ f^{(i)}_\la \right\}$ consists in affecting to   a finite number of cubes $\nu$ the size $e^{(i)}_\la$, and these cubes belong to the dyadic covering of $\la$ so that their scale exceeds the scale of $\la$ by at most $M$. It follows that the Besov norm of $ \left\{f^{(i)}_\la \right\}$ exceeds the one of $\left\{e^{(i)}_\la\right\}$ by, at most, a constant factor. The result for the random sequences follows  because the random variables $\xi^{(i)}_{\la}$ are bounded. 
\end{proof}

The first part of the theorem follows from  the following lemma. 

\BL \label{lem_compact}
Let $A = (a_j)_{j \geq 0}$ be  a positive sequence such that $a_j\rightarrow + \infty$ and $K_A$ be the subset of  $b^{s,q }_p$ defined by the conditions 
$$ \label{lemcompeq}  \cla = 0 \quad \mbox{if} \quad  \la \not\subset [0,1)^d \quad \mbox{ and } \quad  \displaystyle\sum_{j=1}^\infty a_j \left( \sum_{(i, \la) \in \Lambda_j} | \cla 2^{(s-\frac{d}{p})j}|^p \right)^{q/p} \leq 1.    $$
Then $K_A$  is a compact subset of $b^{s,q }_p$. 
\EL

\begin{proof}  
This result is a direct consequence of the compacity of the operator
$$
T : \ell^q \to \ell ^q : (x_j)_{j \geq 0} \mapsto (a^{-1}_j x_j)_{j \geq 0},
$$
see e.g. Chap. 27 of  \cite{MV}. 
\end{proof}
 
%
%

%

%

 
 We will now prove that the  random coefficients \eqref{defCO} are ``large enough and well spread'', which implies maximal divergence for the corresponding wavelet series. 
 We start by proving that  the   $\xi^{(i)}_{j,k}$ cannot be simultaneously small at many successive scales. 
 Let  $j=mM$ be   a multiple of $M$, and for every dyadic cube  $\lambda \subseteq [0,1)^d$ of scale $j$, 
 let $\textgoth{E}_\la$ be the event
\[ \textgoth{E}_\la = \left\{  \forall l \in \{ 1, \dots , L \}, \quad  | \xi^{(i_l)}_{\mu^l (\la) } | \geq \frac{1}{j} \right\} . \] 

\BL \label{lemproprela} The events $\textgoth{E}_\la$ satisfy
 \[ \PP  ( \textgoth{E}_\la) \geq  1-\frac{L}{j} .\]
 If $\textgoth{E}_\la$  holds, then $\forall x \in \la$, there exists $l \in \{1, \dots, L\}$ 
 such that 
 \[ \left| c^{(i_l)}_{\mu^l (\la )} \psi^{(i_l)}_{\mu^l (\la )}  (x) \right| \geq C_0 \frac{e^{(i_l)}_{\la}}{j} . \] 
 Furthermore,  the  $\textgoth{E}_\la$  which correspond to different values of   $j$ (multiples of $M$) are independent.  
 \EL
 
 This lemma follows from  \eqref{denslin}, \eqref{majCo}, \eqref{defCO} and Definition \ref{def_dycov}.
 
 
 \BP \label{prop_div_everywhere}
Almost surely, for every  $x\in [0,1)^d$, the divergence exponent of the wavelet series (\ref{wavexp3}) associated with $\mathcal{C}$ at $x$ satisfies:
$\delta_{ \mathcal C} (x) \in \left[ -s,  -s  + d/p\right]  $,
and, for almost every point $x$ of $[0,1)^d$, 
$\delta_{\mathcal{C}}(x) = -s.$
 \EP
 
 \begin{proof}
 Let us fix a scale  $j$ of the form $j = 2m M$ and  consider a given dyadic cube $\la \subseteq [0,1)^d$ of scale $j$. We denote by $\la_1, \dots  ,\la_m$ the dyadic cubes of scales respectively 
$2mM, (2m-1)M ,\dots  ,(m+1)M$ which contain $\la$. Because of Lemma~\ref{lemproprela}, the probability that none of the $\textgoth{E}_{\la_1}, \dots, \textgoth{E}_{\la_m}$ holds is bounded by 
\[ \frac{L}{2mM}  \frac{L}{(2m-1)M}\cdots \frac{L}{(m+1)M} \leq \left(\frac{L}{(m+1)M} \right)^m \leq e^{-C m \log m} . \] 
Therefore, the probability that at least one of the cubes $\la$ of scale $j= 2mM$ satisfies the property ``{\em none of the $\textgoth{E}_{\la_1}, \dots, \textgoth{E}_{\la_m}$ holds}'' is bounded by $ 2^{dj} e^{-C m \log m}, $
and the Borel-Cantelli Lemma allows to conclude that, almost surely, for $m$ large enough, all cubes $\la$ at scales  $j= 2mM$  satisfy that there is at least one cube $\la'$ of scale $j'$ such that $(m+1)M \leq j' \leq 2mM$ which contains $\la$ and such that $\textgoth{E}_{\la'}$ holds. By Lemma~\ref{lemproprela}, this in turns implies that , for any $ x \in \la$, there exist a cube $\la''$ of scale $j''$ such that $j' < j'' \leq j'+M$ and an index $i$ such that  
 $$ \left| c^{(i)}_{\la'' } \psi^{(i)}_{\la'' }  (x) \right| \geq C_0 \frac{e^{(i)}_{\la'}}{j'} \geq  \frac{C_0}{j'} 2^{-(\log j')^2}2^{-sj'} \geq  \frac{C_0}{j''} 2^{-(\log j'')^2}2^{-sj''}.$$
   It follows that almost surely, $\forall x \in [0,1)^d$ there exists a sequence of cubes $\la_n$ of scales $j_n$ (which grows at most arithmetically) and a sequence of indices $i_n$ such that 
\[ \left| c^{(i_n)}_{\la_n } \psi^{(i_n)}_{\la_n }  (x) \right| \geq \frac{C_0}{j_n} 2^{-(\log j_n)^2}2^{-sj_n}, \] 
so that $\forall x \in [0,1)^d $, $\delta_{\mathcal C} (x) \geq - s$.  This result, together with Proposition \ref{divwase} implies that $\delta_{ \mathcal C} (x) \in \left[ -s,  -s  +d/p \right]  $.

The a.e. divergence rate is $-s$ because,  on one hand,  it is everywhere at least $-s$, and, on other hand, by Proposition \ref{divwase},  $\forall n >0$, the set of points where it is larger than $-s+ 1/n $ has a dimension less that $d-p/n$, hence a vanishing Lebesgue measure. 
 \end{proof}

We will now determine the almost sure wavelet divergence  spectrum of the sequence  ${\mathcal C}$ defined by \eqref{defCO}. We will use  a variant  of Lemma \ref{lemproprela}: Let  $j = mM$ and for every dyadic cube $\lambda \subseteq [0,1)^d$ of scale $j$,  let $\textgoth{F}_\la$ be the event
\[ \textgoth{F}_\la = \left\{  \forall l \in \{ 1, \dots , L  \} ,  \quad  | \xi^{(i_l)}_{\mu^l (\la) } | \geq 2^{-j /\log j}  \right\} . \] 

\BL  The events $\textgoth{F}_\la$ satisfy
 \[  \PP  ( \textgoth{F}_\la) \geq  1-L 2^{-j /\log j}  .\]
 If $\textgoth{F}_\la$  holds, then $\forall x \in \la$, there exists $l \in \{1, \dots, L\}$ 
 such that 
 \[ \left| c^{(i_l)}_{\mu^l (\la )} \psi^{(i_l)}_{\mu^l (\la )}  (x) \right| \geq C_0 e^{(i_l)}_{\la} 2^{-j /\log j}  . \] 
 Furthermore,  the  $\textgoth{F}_\la$  which correspond to different values of   $j$ (multiples of $M$) are independent.  
 \EL 
 
 The proof is similar to the proof  of  Lemma \ref{lemproprela}, using \eqref{denslin}, \eqref{majCo}, \eqref{defCO} and the definition of the dyadic covering property.
We now prove that the divergence exponents of the wavelet series associated to $\mathcal{C}$ is controlled by  the divergence exponents of the deterministic series 
$$
\sum_{i=1}^N \sum_{j \geq 0} \sum_{k \in \ZZ^d} e^{(i)}_{j,k}  \mathbf{1}_{\lambda_{j,k}}. 
$$

\BL \label{lem_E}
Let us denote by $\mathcal{E}$ the sequence whose coefficients are given by (\ref{defCO1}) and by $\mathbf{1}$ the wavelet system define by the functions $\mathbf{1}^{(i)} := \mathbf{1}_{[0,1)^d}$ for every $i \in \{1, \dots, N\}$. 
Almost surely, for every $\gamma$,
$$
\left\{x \in [0,1)^d : \ \delta^{\mathbf{1}}_{\mathcal{E}}(x) \geq \gamma \right\} \subseteq \left\{x \in [0,1)^d: \ \delta^{\Psi}_{\mathcal{C}}(x) \geq \gamma \right\} .
$$
\EL

 \begin{proof}
 Let us fix a scale  $j$ of the form $j = m M$; let $\la \subset [0,1)^d$ be  a given dyadic cube  of scale $j$, and  consider the  sequence of  length $ N_m = [ (\log m)^2 ] +1 $ dyadic cubes of  scales $m M,  (m-1) M, \dots ,(m-N_m+1) M $ that contain $\la$. As done in Proposition~\ref{prop_div_everywhere}, for $m$ large enough, the probability that none of the corresponding $\textgoth{F}_{\la_n}$  holds is bounded by 
\[ \left( L 2^{-m /\log m}  \right)^{(\log m)^2} \leq 2^{-C m \log m}. \] 
A  Borel-Cantelli type argument similar as before yields that almost surely, for $m$ large enough, all cubes $\la$ at scales  $j= mM$  satisfy that there is at least one cube $\la'$ of scale $j' \in \{(m-N_m+1) M , \dots, m M \}$ which contains $\la$ and such that $\textgoth{F}_{\la'}$ holds. Again, it gives a cube  $\la''$ of scale $j'' \in \{j'+1, \dots , j'+M\}$  belonging to the dyadic covering of $\lambda'$ and an index $i$ 
 such that 
 \[ \left| c^{(i)}_{\la''} \psi^{(i)}_{\la''}  (x) \right| \geq C_0 e^{(i)}_{\la'} 2^{-j' /\log j'}  . 
 \] 

Let us now consider $x \in [0,1)^d$ such that $\delta^{\mathbf{1}}_{\mathcal{E}}(x) \geq \gamma$. Then for any $\gamma' < \ga$, there exist a constant $C>0$ and infinitely many indices $(j_n,k_n)$ for which
$$
x \in \lambda_{j_n,k_n} \quad \mbox{ and } \quad \left| e^{(i)}_{j_n,k_n}\right| \geq C 2^{ \ga' j_n} \quad \forall i \in \{1, \dots, N\}.
$$
For every $n$, let $m_n$ be the integer satisfying $m_n M \leq j_n < (m_n+1)M$ and $\lambda_n$ be the dyadic cube of scale $m_nM$ which contains $\lambda_{j_n,k_n}$. Then, if $n$ is large enough, we have obtained that there is a cube $\lambda'_n$ of scale $j'_n \in \{(m_n-N_{m_n}+1) M , \dots, m_n M \}$ which contains $\la_n$, a cube  $\la''_n$ of scale $j''_n \in \{j'_n+1, \dots , j'_n+M\}$  belonging to the dyadic covering of $\lambda'$ and an index $i_n$  such that 
 \[ \left| c^{(i_n)}_{\la''_n} \psi^{(i_n)}_{\la''_n}  (x) \right| \geq C_0 e^{(i_n)}_{\la'_n} 2^{-j'_n /\log j'_n}  . 
 \] 
Thanks to the hierarchical property of the sequence $\mathcal{E}$ given in Remark \ref{hierar}, up to a factor $2^{\beta(j_n - j'_n)}$ which does not affect the divergence exponent since the difference $j_n-j'_n$ is of order $(\log (j_n))^2$, the size of $e^{(i_n)}_{\la'_n}$ is comparable with the size of $e^{(i_n)}_{j_n,k_n}$. This gives then a control of the size of $| c^{(i_n)}_{\la''_n} \psi^{(i_n)}_{\la''_n}  (x) |$ and since $j''_n - j'_n \leq M$, we get that $\delta^{\Psi}_{\mathcal{C}}(x) \geq \gamma'$. The conclusion follows.
 \end{proof}

\BP \label{prop_spectrum}
Almost surely, for any non-empty open set  $B\subset [0,1)^d$,
 \BE  
\dim\left( \{x \in B : \delta^{\Psi}_{\mathcal{C}}(x) = \gamma\}\right)  =  d -sp -\ga p  .  \EE
\EP

\begin{proof}
Let $\al \in  [1, +\infty ) $ and denote by  $E_\al$ the set of points $x$ of $[0,1)^d$ for which there exists infinitely many indices $(j_n,k_n)$ with $j_n \to + \infty$ such that
$$
0 \leq  x - \frac{k_n}{2^{\lfloor j_n/\al \rfloor}} < \frac{1}{2^{j_n}}.
$$
We denote by $\la_n$ the  sequence of dyadic cubes of scale $j_n$ that contain $x$.  If $x\in E_\al$, then for any index $i$,  
\[ e^{(i)}_{\la_n }  \geq 2^{\left(\frac{d}{p} -s\right)j_n}  2^{-\frac{d}{p\al } j_n} 2^{-\frac{j_n}{\log j_n } } , \]
so that 
\[ \delta^{\mathbf{1}}_{\mathcal E} (x) \geq \gamma (\al) := \frac{d}{p } -s -\frac{d}{p\al }.\]
Lemma \ref{lem_E} implies then that
$$
E_\al \subseteq \left\{x \in [0,1)^d: \ \delta^{\Psi}_{\mathcal{C}}(x) \geq \gamma (\al) \right\}.
$$
The computation of the Hausdorff dimension of $E_\al$ is standard; we will need  the   definition of the { \em modified Hausdorff measure}. If   $A\subset\RR^d$, 
  $\varepsilon>0$  and  $\delta \in \RR $,  let 
\[ M^{\delta, \gamma}_{\varepsilon} (A) =\inf_R \;  \left( \sum_{  i} | A_i |^\delta | \log ( | A_i |)|^\gamma   \right) ,\]
where $R$ denotes an  $\ep$-covering of   $A$, i.e. a covering of  $A$  by bounded sets  $\{ A_i\}_{i \in \NN}$
of diameters   $| A_i | \leq \varepsilon$.
The infimum  is therefore taken on all possible  $\ep$-coverings.
For any  $\delta \in \ ]0,d]$ and $\gamma \in \RR$, the quantity 
\[ \mathcal{H}_{\delta,\gamma} (A) = \displaystyle\lim_{\ep\rightarrow 0^+}  M^{\delta, \gamma}_{\ep} (A)  \]
defines  the $(\delta, \gamma )$-dimensional (outer) Hausdorff measure of 
$A$.
A standard result of dyadic approximation 
(for instance, see the mass
transference principle of \cite{BV06} and Theorem 2 of \cite{J00}),  
yields  that, for any ball $B$ of  non-empty interior,  
 \[ \mathcal{H}_{d/\al, 2} ({ E}_{\al}  \cap B ) >0  .\]
  On other hand,  
 (\ref{maxdiv}) implies that the set of points for which the wavelet series associated to $\mathcal{C}$ has a divergence exponent strictly  larger than  $\ga (\al) $ have a $(d/\al, 2) $-Hausdorff measure that vanishes (it follows from interpreting it as a countable union of the set of points where the divergence exponent is   larger than  $\ga (\al) + 1/n$, which, by Proposition \ref{divwase} have a Hausdorff dimension bounded by $\frac{d}{\al} -\frac{p}{n}$ and therefore   have a vanishing $(d/\al, 2)$-Hausdorff measure). 
 It follows that the set of points with a divergence exponent    $\ga (\al) $ has  positive $(d/\al, 2) $-Hausdorff measure, and thus has a Hausdorff dimension equal to $d/\al$. 
 \end{proof}

\begin{proof}[\bf Proof of Theorem \ref{thm1}]
It is a direct consequence of Lemma \ref{lem_compact}, Proposition \ref{prop_div_everywhere} and Proposition \ref{prop_spectrum}.
\end{proof}

The next remark  will be used in the proof of Theorem \ref{th1}.

\BR \label{rem_max}
A careful inspection of the previous proofs shows in particular that if $\mathcal{D} = \left\{ d^{(i)}_\la \right\}$ is a sequence such that, for every $x \in [0,1)^d$, there are infinitely many cubes $\lambda_n$ containing $x$ and of increasing scales $j_n$, and infinitely many indices $i_n$ such that
$$\frac{j_{n}-j_{n-1}}{j_n} \to 0 \qquad \mbox{and} \qquad \left| d^{(i_n)}_{\lambda_n} \psi^{(i_n)}_{\lambda_n}(x) \right| \geq C e^{(i_n)}_{\lambda_n},$$
then the wavelet series associated to $\mathcal{D}$ is maximally divergent in $[0,1)^d$.
\ER

\BR \label{rem_tilde} We have constructed a random sequence  ${\mathcal C}$ whose sample paths almost surely yield a maximally divergent wavelet series on $[0,1)^d$.  In order to obtain the same result on the whole $\RR^d$, it  suffices to consider the process
 \[ \widetilde{\mathcal C} = \sum_{k \in \ZZ^d} e^{-k} \; {\mathcal C}^k  , \] 
where the ${\mathcal C}^k$ are independent copies of ${\mathcal C}$ translated by $k$. 
\ER

We end this section by the following corollary, which states that the results obtained for the sequence ${\mathcal C}$ also hold for ``translates''  of ${\mathcal C}$.  This result will be the key step to obtain generic results of maximal divergence in the sense of prevalence below.

\BC \label{estidecentered} 
Let ${ \mathcal D}$ be a given sequence in 
$b^{s,q }_p$. Then  the wavelet series associated to ${\mathcal C} + { \mathcal D}$   is almost surely  maximally divergent in $[0,1)^d$.
\EC

\begin{proof}
This corollary  is a direct consequence of the fact that all the properties that we have obtained for  ${\mathcal C} $ also hold for  ${\mathcal C}  + { \mathcal D}$. Indeed the only difference is that the random variables defining its coefficients are no more centered, but are shifted by  deterministic quantities. Beside their independence (which still holds),  the only  property that we used for these coefficients is  lower bounds for the quantities $ \PP ( \{Ê| \cli | \leq  c\} )$, but these quantities become smaller when the densities are shifted, as a straightforward consequence of the explicit density  \eqref{denslin}  that we have chosen  (because this density is  an even function, which is non-increasing on $\RR^+$). 
\end{proof}


\subsection{Prevalence} 

The notion of prevalence supplies an extension of the notion of ``almost everywhere'' (for the  Lebesgue measure)  in an infinite dimensional setting. In
a metric infinite dimensional vector space,  no measure  is both  $\sigma$-finite  and translation invariant. However, a natural notion of ``almost everywhere'' which is translation invariant  can be defined as follows,  see  \cite{Christensen74,HSY92}; ``zero-measure sets'' 
thus defined  are called Haar-null.

\BD\label{def:prev}
Let E be a complete metric vector space. A Borel set $A\subset E$ is  {Haar-null} if there exists a compactly supported probability measure $\mu$  
such that
\begin{equation}\label{trans}
\forall x\in E,\quad \mu(x+A)=0. 
\end{equation} 
  If this property holds, the measure $\mu$ is said to be {transverse} to $A$.

A subset of E is called {Haar-null} if it is contained in a Haar-null Borel set. The complement of a Haar-null set is called a {prevalent set}.
\ED

The following results of \cite{Christensen74} and  \cite{HSY92}  enumerate  important properties of prevalence:
\begin{itemize}
\item If S is Haar-nul, then  $\forall x\in E$, $x+S$ is Haar-nul.
\item If the dimension of $E$ is finite, S is Haar-null if and only if $S$ has  Lebesgue measure $0$.
\item Prevalent sets are dense.
\item The intersection of a countable collection of prevalent sets is prevalent.
\item If the dimension of $E$ is infinite,  compact subsets of  $E$ are Haar-nul.
\end{itemize}

\BR In order to prove that a set is Haar-nul, one can often use for transverse measure the Lebesgue measure  on the unit ball
of  a finite dimensional subset V; Condition (\ref{trans}) becomes 
\[ \forall x\in E, \;\;\;\; (x+A)\cap V \;\;\mbox{  is of Lebesgue measure
zero}.\]
 In this case $V$ is called a {probe} for the complement of $A$.
\ER

\BR If $E$ is a function space,  choosing a probability measure on $E$ is equivalent to choosing a random process $X_t$ the
sample paths  of which  almost surely belong to a compact subset of $E$. Thus,  the definition of a Haar-null set can be rewritten as  follows:
Let
$\mathcal{P}$ be a property that can be satisfied by points of  $E$ and let  
\[ A= \{f\in E:\; \mathcal{P}(f) \;
\mbox{holds } \}. \]
 The condition $\mu(f+A)=0$ means the event $\mathcal{P}(X_t-f)$ has  probability zero. Therefore, a way
 to check that  a property
$\mathcal{P}$ holds  only on a Haar-null set  is to  exhibit  a random process $X_t$ whose  sample paths  are in a compact subset of E
and is such that 
\begin{displaymath}
\forall f\in E,\; a.s. \; \mathcal{P}(X_t+f)\;\mbox{does not hold.} 
\end{displaymath}   
\ER

With a slight abuse of language, when a property holds on a prevalent set, we will say that it holds almost everywhere. 
In the following theorem, we obtain the  optimality of Proposition \ref{divwase} in a prevalent sense.

\BT \label{th1prev} Let $s\in \RR$ and $p, q \in (0, + \infty]$. If the wavelet system $\left\{\pjk\right\}$ satisfies the dyadic covering property, then the wavelet series of almost every sequence of $b^{s,q}_p$ is maximally divergent in $\RR^d$. 
\ET 

\begin{proof}
Let us consider the random sequence ${\mathcal{C}}$ introduced in (\ref{defCO}). Corollary \ref{cor1} yields  that for every given  sequence  $\mathcal{D}$ of $b^{s,q}_p$, almost surely, the wavelet series associated to $\mathcal{D} + {\mathcal{C}} $ is maximally divergent in $[0,1)^d$. Moreover, it follows from Theorem \ref{thm1}, that the random  sequence ${\mathcal C}$ takes values in a compact subset of $b^{s,q }_p$. Hence the result when we consider the divergence at points of $[0,1)^d$. Since a countable intersection of prevalent sets is prevalent, the conclusion is obtained by covering $\R^d$ with a countable family of cubes.  
\end{proof}

\section{Other generic results of divergence} 
\label{sec:generic}

As soon as an object with some extremal behavior has been exhibited, it is natural to wonder if this  behavior  is generic, in some sense. Prevalence allows to show  that the set formed by these special objects is large in a measure sense. One  can also consider the Baire categories genericity, called residuality: There exists a dense $G_\delta$ set of elements sharing this behavior in a well chosen topological vector space.  More recently, the algebraic structure of this kind of sets has also been investigated, using the notion of lineability. 

In this section, we obtain the equivalent of Theorem \ref{th1prev} in the senses supplied by Baire categories and lineability.

\subsection{Residuality}

The aim of this section is to construct a countable intersection of dense open sets of $b^{s, q}_p$  the elements of which have maximally divergent wavelet series. 

\BT \label{th_res}
Let $s \in \RR$ and  $p,q \in (0, + \infty]$. If the wavelet system $\left\{\pjk\right\}$ satisfies the dyadic covering property, then the set of sequences of $b^{s, q}_p$ whose  wavelet series is  maximally divergent in $\RR^d$ is residual in  $b^{s, q}_p$.
\ET

\begin{proof}

As usual, it is sufficient to  consider  divergence rates at points of the unit cube.
 The idea of the construction is that any sequence whose  coefficients are sufficiently close to those of the sequence ${\mathcal E}$, whose wavelet coefficients are defined by (\ref{defCO1}), on infinitely many dyadic cubes will have a maximally divergent wavelet series. Let us study separately different cases.

\medskip

\underline{The case $p,q < + \infty$ :} 

\medskip

We consider the norm
\[
 \| \mathcal{D}\| = \left( \sum_{j \geq 0} \Bigg(\sum_{(i,\lambda) \in \Lambda_j} |d^{(i)}_{\lambda}2^{(s-\frac{d}{p})j}|^p \Bigg)^{q/p} \right)^{1/q} 
\] 
on $b^{s,q}_p$. Since $p,q < + \infty$, the space  $b^{s, q}_p$ is separable and finite sequences with rational coefficients form a dense subspace in this space. Let $(\mathcal{F}_n)_{n \in \NN}$ denote such a dense  sequence.  For every $n \in \NN$, there is $N_n$ such that the coefficients of $\mathcal{F}_n$ are equal to $0$ at scales $j \geq N_n$. Without loss of generality, one can assume that the sequence $(N_n)_{n \in \NN}$ is increasing. For every $n \in \NN$, let us define $$\mathcal{G}_n = \mathcal{F}_n + \frac{1}{N_n}{\mathcal E} \, .$$
 By construction, the sequence $(\mathcal{G}_n)_{n \in \NN}$ is dense in $b^{s,q}_p$.  Finally, let us consider the set
\[
\mathfrak{R} = \bigcap_{m \in \NN} \bigcup_{n \geq m} B(\mathcal{G}_n,r_n),
\]
where $$r_n = \frac{1}{2N_n} 2^{-(\log (N_n+M))^2}2^{-\frac{d (N_n +M)}{p}}.$$
The set $\mathfrak{R}$ is clearly a countable intersection of dense open sets. Moreover, let us remark that if $\left\{d^{(i)}_{\lambda}\right\}$ denotes the  coefficients of a sequence $\mathcal{D}$ in $B(\mathcal{G}_n,r_n)$, then 
for every $(i,\lambda) \in \Lambda_j$ with $j \geq N_n$, 
\begin{equation}\label{eq_boule}
 \big|d^{(i)}_{\lambda}- \frac{1}{N_n}e^{(i)}_{\lambda} \big| < 2^{(\frac{d}{p}-s)j} r_n.
\end{equation}

If $\mathcal{D}$ belongs to the residual set $\mathfrak{R}$, it belongs to infinitely many balls $B(\mathcal{G}_{n_l},r_{n_l})$. Clearly, the dyadic covering property and inequality (\ref{eq_boule}) give that at every $x$, the divergence exponent is at least equal to $-s$, and as done at the end of the proof of Proposition \ref{prop_div_everywhere}, that it is equal to $-s$ for almost every $x$. Let us now compute the divergence spectrum of $\mathcal{D}$. For $\alpha \geq 1$, let us denote by $E_{\alpha}(\mathcal{D})$ the set of points $x$ satisfying the following property: for infinitely many $l$, there exists $k_{n_l}$ such that
$$
0 \leq  x - \frac{k_{n_l}}{2^{\lfloor N_{n_l}/\al \rfloor}} < \frac{1}{2^{N_{n_l}}}.
$$
Remark that the  definition of $E_{\alpha}(\mathcal{D}) $ is very similar to that of $E_{\alpha}$ introduced in Proposition \ref{prop_spectrum}, except that only a subsequence of the integers is used. Consequently, if $x \in E_{\alpha}(\mathcal{D})$, we have
\BE \label{estim} 
e^{(i)}_{\la_{n_l} }  \geq 2^{\left(\frac{d}{p} -s\right)N_{n_l}}  2^{-\frac{d}{p\al } N_{n_l}} 2^{-\frac{N_{n_l}}{\log N_{n_l} } } , \EE
for any index $i$, where  $\la_{n_l}$ denote the dyadic cubes of scale $N_{n_l}$ that contain $x$. Using the dyadic covering property, there exists a dyadic cube $\mu$ of scale $j \in \{N_{n_l} +1, \dots, N_{n_l} + M\}$ and an index $i$ such that $|\psi^{(i)}_{\mu}(x)| \geq C_0$. Together with (\ref{eq_boule}), it gives
\begin{eqnarray*}
|d^{(i)}_{\mu}\psi^{(i)}_{\mu}(x)| & \geq & \frac{1}{N_{n_l}} e^{(i)}_{\mu} -  2^{(\frac{d}{p}-s)j}\frac{1}{2N_{n_l}}2^{-(\log (N_{n_l}+M))^2}2^{-\frac{d (N_{n_l}+M) }{p}} \\
& \geq & C_0 \frac{1}{2j} e^{(i)}_{\mu}.
\end{eqnarray*}
Remark \ref{hierar} concerning the hierarchic property of the coefficients $\left\{ e^{(i)}_{\la} \right\}$ and inequality (\ref{estim}) imply that  the divergence rate of $\mathcal{D}$ at $x$ is larger than $\frac{d}{p} - s - \frac{d}{\alpha p}$. The computation of the Hausdorff dimension of $E_{\alpha}(\mathcal{D})$ is then computed as done for $E_{\alpha}$ in Proposition~\ref{prop_spectrum}. It follows that the set of points where the divergence rate of $\mathcal{D}$ is exactly $\frac{d}{p} - s - \frac{d}{\alpha p}$ has dimension $\frac{d}{\alpha}$. The theorem follows. 

\medskip

\underline{The case $p,q = + \infty$ :}

\medskip

 This case corresponds to the H\"older space $C^s$. We consider the norm
\[
 \| \mathcal{D}\| = \sup_{j \geq 0} \sup_{(i,\lambda) \in \Lambda_j} 2^{sj}|d^{(i)}_{\lambda}|
\] 
on $C^s$. For each integer $n$, let us denote by $E_n$ the set of sequences of $C^s$ whose coefficients are each a non-vanishing multiple of $2^{-sj-n}$. If $\mathcal{D} \in E_n$, then its coefficients $\left\{d^{(i)}_{\lambda} \right\}$ satisfy
\[
|d^{(i)}_{\lambda}| \geq 2^{-n-js}.
\]
Let us set
\[
A_n = E_n + B(0, 2^{-n-1}) 
\]
and let us notice that the coefficients $\left\{a^{(i)}_{\lambda}\right\}$ of any sequence $\mathcal{A} \in A_n$ satisfy
\[
 |a^{(i)}_{\lambda}| \geq 2^{-n-sj-1} .
\]
Using this last inequality, Proposition \ref{maxdiv} and the dyadic covering property, one directly obtains that the divergence exponent of the wavelet series associated with $\mathcal{A}$ is exactly $-s$ at every point. 

Let us now define
\[
\mathfrak{R} = \bigcup_{n \in \N} A_n.
\]
If $\mathcal{B} \in C^s$ has  coefficients $\left\{b^{(i)}_{\lambda}\right\}$ and if we define for every $n \in \N$, the sequence $\mathcal{A} \in A_n$  by
\[
a^{(i)}_{\lambda}  = \left\{
\begin{array}{ll}
2^{-sj-n} \lfloor 2^{sj+n} b^{(i)}_{\lambda} \rfloor & \text{ if } \lfloor 2^{sj+n} b^{(i)}_{\lambda} \rfloor \neq 0, \\
2^{-sj -n} & \text{ otherwise,} 
\end{array}
\right.
\]
then $\| \mathcal{B}- \mathcal{A} \| \leq 2^{-n}$. Therefore, $\mathfrak{R}$ is a dense open set (hence a $G_\delta$ set) of $C^s$ whose element have a divergence exponent equal to $-s$ at each point, and the conclusion follows.

\medskip

\underline{The cases $p= + \infty, q < + \infty$ and $p < + \infty, q = + \infty$ :}

\medskip
It suffices to proceed as in the previous cases, with obvious adaptations.






\end{proof}



\subsection{Lineability}

In a nutshell, proving a generic result in the sense of lineability consists in proving that this result holds for {\em every} (non zero) element of a subspace of infinite dimension (the cardinality of the space allowing for different variants of the notion). We will define this subspace by defining explicitly an uncountable Hamel basis (indeed, the space considered in the definition of lineability is understood as  a space generated by the {\em finite} linear combinations of the elements of the basis).

\BD
Let $X$ be a vector space, $M$ a subset of $X$, and $\kappa$ a cardinal number. The subset $M$ is said to be {$\kappa$-lineable} if $M \cup \{0\}$ contains a vector subspace of dimension \(\kappa\). The set $M$ is {lineable} if the existing subspace is infinite dimensional. When $X$ is a topological vector space and when the above vector space can be chosen to be dense in $X$, we  say that $M$ is {$\kappa$-dense-lineable} (or, simply, dense-lineable if $\kappa$ is infinite).
\ED

This recent concept  has attracted the attention of many authors, see e.g. the detailed  review of L. Bernal-Gonz\'alez, D. Pellegrino and J.B. Seoane-Sep\'{u}lveda \cite{BPS}.
Recently, L. Bernal-Gonz\'alez \cite{BG10} introduced the notion of {\em maximal lineability} (and that of {\it maximal dense-lineability}) meaning that the dimension of the existing vector space is equal to the dimension of $X$.

The aim of this subsection is to obtain the equivalent of Theorem~\ref{th_res} and Theorem~\ref{th1prev} in the context of lineability. First, let us introduce the sequences which will form the basis of the linear subspace. For every $a>0$, let us consider the sequence $\mathcal{E}_a = \left\{ (e_{a})^{(i)}_{\lambda}\right\}$, with 
$$
(e_{a})^{(i)}_{\lambda} = \frac{1}{j^a} e^{(i)}_{\lambda}
$$
for every $(i,\lambda) \in \Lambda_j$, where the coefficients $\left\{ e^{(i)}_{\lambda}\right\}$ are defined by (\ref{defCO1}). Since the sequence $\left\{ e^{(i)}_{\lambda}\right\}$ belongs to $b^{s,q}_p$, it is clear that $\mathcal{E}_a \in b^{s,q}_p$ for every $a>0$.  
The result of lineability will direcly follow from the next straightforward lemma which states that the coefficients of any linear combination of the sequences $\mathcal{E}_a$ are of the order of magnitude of its ``largest'' component. Therefore, it has the same divergence properties as this component. Moreover, it gives the linear independence of the sequences $(\mathcal{E}_a)_{a >0}$ since no non-zero linear combination of these sequences can be identically equal to $0$. 

\BL \label{lem1} Let $s\in \RR$ and $p, q \in (0, + \infty]$. Let $ n \geq 1,$ $a_n> \dots > a_1   > 0 $, and  $ k_1 , \dots, k_n  \neq 0 $, and consider the sequence 
$$\label{fung} \mathcal{D} = \sum_{i = 1}^n k_i  \mathcal{E}_{a_i}.  $$
Then, the coefficients of $\mathcal{D}$ satisfy that 
$$  \label{propCO}  
\frac{| k_1| }{2 j^{a_1}} e^{(i)}_{\lambda} \leq |  d^{(i)}_{\lambda}|  \leq   \frac{2 | k_1 |}{j^{a_1}} e^{(i)}_{\lambda}  $$
for every $(i,\lambda) \in \Lambda_j$, with $j$ large enough. 
\EL

\BT \label{th1}
Let $s \in \RR$ and  $p,q \in (0, + \infty]$. If the wavelet system $\left\{\pjk\right\}$ satisfies the dyadic covering property, then the set of sequence of $b^{s, q}_p$ whose  wavelet series is  maximally divergent in $\RR^d$ is maximal lineable in $b^{s, q}_p$.
\ET

\begin{proof}
Using the dyadic covering property, Remark \ref{rem_max} and Lemma \ref{lem1}, one directly gets that any non-zero sequence in the linear span of the sequences $\mathcal{E}_a$ has a maximally divergent wavelet series. 
\end{proof}


Actually, in the separable case, one can slightly modify the above construction in order to get a dense subspace of sequences with maximally divergent wavelet series. It is given by the next result.

\BC \label{cor1}
Let $s \in \RR$ and  $p,q \in (0, + \infty)$. If the wavelet system $\left\{\pjk\right\}$ satisfies the dyadic covering property, then the set of sequence of $b^{s, q}_p$ whose  wavelet series is  maximally divergent is maximal dense-lineable in $b^{s, q}_p$.
\EC

\begin{proof}
As done in the proof of Theorem \ref{th_res}, let $(\mathcal{F}_n)_{n \in \NN}$ denote the dense sequence of finite sequences with rational coefficients. Let us also choose a sequence $(a_n)_{n \in \NN}$ of different positive numbers. For every $n \in \NN$, we fix $\epsilon_n>0$ such that
\[
\| \epsilon_n \mathcal{E}_{a_n} \|_{b^{s,q}_p} < \frac{1}{n},
\]
and we define $\mathcal{G}_n = \mathcal{F}_n + \epsilon_n \mathcal{E}_{a_n} $. By construction, the sequences $\mathcal{G}_n$, $n \in \NN$, form a dense subspace of $b^{s, q}_p$.  Finally, we consider the subspace $\mathfrak{D}$ generated by the sequences $\mathcal{G}_n$, $n \in \NN$, and the sequences $\mathcal{E}_a$, $a \in A$, where $A = \{ a > 0 : a \neq a_n \ \forall n \in \NN\}$. Clearly, since it contains the sequences $\mathcal{G}_n$, $n \in \NN$, the subspace $\mathfrak{D}$ is dense in $b^{s, q}_p$. Moreover, it has maximal dimension since it contains the linearly independent sequences $\mathcal{E}_a$, $a \in A$. It remains to prove that any non-zero element of $\mathfrak{D}$ has maximally divergent wavelet series. Such an element has, for large scales, the same  coefficients as a non-zero linear combination of the sequences $\mathcal{E}_a$, $a>0$. The conclusion follows then with the same arguments as in Theorem \ref{th1}.
\end{proof}



\BR  In contradistinction with the Baire or prevalence case,  the notion of lineability is not stable under intersection. Therefore different lineability results can hold simultaneously: This is the case  for Theorem \ref{th1}. In order to be  clearer we now mark the dependence in $s$ and $p$ of the sequence $\mathcal{E}$, whose coefficients are given by  (\ref{defCO1}), and  denote it  by $\mathcal{E}^{s,p}$.  We  now pick $s', p'$ such that $s'-d/p' > s-d/p$ and $p \geq p'$. Classical Besov embeddings give $b^{s',q}_{p'}\subseteq b^{s,q}_p$. Then, the subspace constructed in Theorem \ref{th1} or in Corollary \ref{cor1} can also be constructed starting with the sequence $\mathcal{E}^{s',p'}$ instead of $\mathcal{E}^{s,p}$, and (\ref{dimgenerloc}) is replaced by $ \dim \left(  \{x \in B : \delta_{\mathcal{C}}(x) = \gamma\}\right)  =  d -s'p' -\ga p'$. 
\ER

\section{Maximal divergence at a given point}\label{sec:point}

We  now determine the generic divergence exponent at a given (fixed) point of wavelet series associated with sequences in a Besov space.

\BP \label{divmaxpoint}
 Let $s\in \RR$,  $p, q \in (0, + \infty]$ and $x_0\in \RR^d$. 
There exist a sequence  ${\mathcal C} \in b^{s,q}_p $  whose divergence exponent at $x_0$ satisfies 
\[ \delta_{{\mathcal C}} (x_0) = -s + \frac{d}{p}. \] 
\EP

\begin{proof} Let $x_0 \in \RR^d$ be fixed. Using the dyadic covering property, let us consider a sequence of cubes $\lambda_n$ of strictly increasing scales and a sequence of indices $i_n$ such that
\[
\left|\psi_{\lambda_n}^{(i_n)}(x_0) \right| \geq C_0.
\]
Let us define the sequence $\mathcal{C} = \left\{ \cla \right\}$ by setting  for every $(i,\lambda) \in \Lambda_{j}$  
\BE \label{defga}
\cla = \left\{
\begin{array}{ll}
 2^{-(\log  j)^2}  2^{-(s-d/p)j}  & \text{ if }  \exists n \ : (i,\lambda)=(i_n,\lambda_n) \\
 0 & \text{ otherwise}.
\end{array}
\right.  \EE 
Since there is at most one non-vanishing wavelet coefficients at each scale, 
\[ \left( \sum_{(i, \lambda) \in \Lambda_j} | \cla 2^{(s-\frac{d}{p})j}|^p \right)^{1/p}\leq   2^{-(\log  j)^2} 
\]
and ${\mathcal C}$ clearly belongs to the Besov space $b^{s,q}_p$. Furthermore, the choice of $(i_n,\lambda_n)$ implies that the divergence exponent of $\mathcal{C}$ at $x_0$ is larger than $-s + \frac{d}{p}$, and Proposition \ref{divwase} gives the conclusion. 
\end{proof}

Let us now extend Proposition \ref{divmaxpoint} in the setting of genericity, using the three different notions.

\BT \label{divmaxpointres}
Let $s\in \RR$,  $p, q \in (0, + \infty]$ and $x_0\in \RR^d$. The set of sequences $\mathcal{D}\in b^{s,q}_p$  such that  
\[ \delta_{\mathcal{D}} (x_0) = -s + d/p \] 
is prevalent,  residual and maximal lineable in $b^{s, q}_p$.
\ET

\begin{proof}
Let us start with the notion of prevalence. We will use for probe the space generated by the sequence $\mathcal{C}$ whose coefficients are given by (\ref{defga}) and prove that, $\forall \be < -s+d/p $, and  for any sequence $ \mathcal{D}\in b^{s,q}_p$, the line $ (\mathcal{D} + a \mathcal{C})_{a \in \RR}$ contains at most one sequence whose wavelet series diverges at rate at most $\be$ at $x_0$.

Indeed, suppose that there exist two such sequences corresponding to two values $a_1$ and $a_2$. Then, by the triangular inequality and the definition of $\mathcal{C}$, there are infinitely many scales $j$ for which there is $(i, \lambda) \in \Lambda_j$ such that
\[   | a_1-a_2 | 2^{-(\log j)^2} 2^{-(s-d/p)j} | \pla (x_0)  | \leq  C 2^{\be j} \]
and $|\pla(x_0)| \geq C_0 >0$.
Consequently, for infinitely many $j$,
\[ C_0  | a_1-a_2 | \leq  C  2^{(\be + s-d/p)j} 2^{(\log j)^2} \]
and making $  j \rightarrow + \infty$, we get that $a_1 = a_2$.

\medskip
The proof of  residuality relies on arguments similar to those used in the proof of Theorem \ref{th_res}: in the case $p,q< + \infty$, it suffices to construct the sequence $(\mathcal{G}_n)_{n\in \N}$ via the coefficients (\ref{defga}) instead of (\ref{defCO1}) and to consider as radius $r_n = \frac{1}{2N_n} 2^{-(\log(N_n+M))^2}$. Let us notice that one can assume that the sequence of scales $j_n$ of the cubes which appear in the construction (\ref{defga}) grows at most arithmetically (with common difference $M$).  As done previously, it ensures that the wavelet coefficients of the elements of $B(\mathcal{G}_n,r_n)$ are closed enough to those defined by (\ref{defga}). 

 The other cases follow  from straightforward  modifications  of these arguments.

\medskip
As regards lineability, consider for every $a>0$, the sequence $\mathcal{C}_a = \left\{ (c_{a})^{(i)}_{\lambda}\right\}$, with 
$$
(c_{a})^{(i)}_{\lambda} = \frac{1}{j^a} c^{(i)}_{\lambda}
$$
for every $(i,\lambda) \in \Lambda_j$, where the coefficients $\left\{ c^{(i)}_{\lambda}\right\}$ are defined by (\ref{defga}). As in Lemma~\ref{lem1}, we notice that the coefficients of any non zero linear combination of the $\mathcal{C}_a$, $a>0$, are of the order of magnitude of $\frac{1}{j^{a_1}}\cla$ for some $a_1>0$. The conclusion follows then directly from Proposition~\ref{divmaxpoint}.
\end{proof}

 \bibliographystyle{plain}
 \bibliography{Lineab-1}

\end{document}